\documentclass{amsart}
\usepackage{amssymb,latexsym} 

\begin{document}
\title{\bf On functions without a normal order}
\author{Peter Shiu}
\address{353 Fulwood Road, Sheffield, S10 3BQ, United Kingdom}
\email{p.shiu@yahoo.co.uk}
\date{13 June 2016}
\begin{abstract} 
The method of Tur\'an in establishing the normal order for the number of prime divisors of a number is used to show that a certain class of arithmetic functions do not have a normal order. \end{abstract}
\maketitle

\def\cc#1{\mathcal{C}_#1}
\def\dd#1{\mathcal{D}_#1}
\def\ee#1{\mathcal{E}_#1}
\def\jj#1{\mathcal{J}_#1}
\def\xx#1{\mathcal{X}_#1}
\def\qq#1{\mathcal{Q}_#1}
\def\mm{\mathcal{M}}
\def\nd#1#2{\hbox{$#1\!\!\not|#2$}}
\def\nds#1#2{\hbox{$\scriptscriptstyle#1\not\kern1.5pt|#2$}}
\section{Introduction}
The normal order of an arithmetic function, defined in \cite[p.~356]{HandW}, 
measures the `usual~size' of the function: 
A function $\psi(n)\ge0$ is said to have a normal order $f(n)$ if, to every $\epsilon>0$, the number of $n\le x$ for which $|\psi(n)-f(n)|<\epsilon f(n)$ is~$o(x)$, as $x\to\infty$. 
It is tacitly assumed that $f(n)$ is increasing---otherwise, every such $\psi(n)$ has itself as normal order. 

The notion was first introduced by G.~H.~Hardy and S.~Ramanujan~\cite{HandR}, who proved that $\omega(n)$, the number of distinct prime divisors of~$n$, has the normal order~$\log\log n$. 
Their proof was much simplified by P.~Tur\'an (\cite[p.~356]{HandW}, \cite{pt}), who showed that the result can be established from the asymptotic formulae for the first and the second moments of~$\omega(n)$; indeed it is sometimes said that probabilistic number theory stems from~\cite{pt}. 
By applying Tur\'an's method `in reverse', so to speak, S.~L.~Segal~\cite{sls} showed that Euler's totient function~$\phi(n)$ does not have a normal order. 
We distil the argument used by~Segal, thereby extending his result to a certain class of arithmetic functions.

\section{A class of functions without a normal order}
Let $\mm$ denote the class of arithmetic functions $\psi$ for which there are positive constants $A,B,C$ such that $0\le\psi(n)<Cn$ and, as $x\to\infty$, 
\begin{equation}\label{moments}
\sum_{n\le x}\psi(n)\sim\frac{Ax^2}2\qquad\text{and}\qquad\sum_{n\le x}\psi^2(n)\sim\frac{Bx^3}3.\end{equation} 

{\bf Theorem}. \emph{Let $\psi\in\mm$. If $A^2<B$ then $\psi$ does not have a normal order. }

\medskip
\emph{Proof}. 
Let $A,B,C$ be constants associated with $\psi\in\mm$, and set
\begin{equation}\label{rx}
R(x)=\sum_{n\le x}\Big(\psi(n)-An\Big)=o(x^2), \qquad\text{as}\quad x\to\infty.
\end{equation}    
Suppose that $\psi(n)$ has the normal order~$f(n)$; we may assume without loss that $f(n)<2Cn$, so that $|\psi(n)-f(n)|\le\max\{\psi(n),f(n)\}<2Cn$. 
Making use of (\ref{rx}), and $f(n)$ being increasing, we find, by partial summation, that 
\begin{align}\label{partial}
\Big|\sum_{n\le x}\Big(\psi(n)-An\Big)f(n)\Big|
&\le\max_{n\le x}|R(n)|\Big\{\sum_{n\le x-1}\Big(f(n+1)-f(n)\Big)+f(x)\Big\}\\
&=o(x^3)\quad\text{as}\quad x\to\infty. \notag
\end{align}

Let $\epsilon>0$. 
Appealing to the definition of normal order and separating terms 
depending on whether $|\psi(n)-f(n)|<\epsilon f(n)$, or not, we then have, as $x\to\infty$, 
\begin{equation}\label{variance}
\sum_{n\le x}(\psi(n)-f(n))^2\le 4\epsilon^2C^2\sum_{n\le x}n^2+4C^2x^2o(x)
=\frac{4\epsilon^2C^2x^3}3+o(x^3).\end{equation}
From (\ref{moments}), (\ref{partial}), (\ref{variance}), together with 
\begin{align*}
\psi^2(n)&=A^2n^2+(\psi(n)-f(n))^2+2(\psi(n)-An)f(n)-(f(n)-An)^2\\
&\le A^2n^2+(\psi(n)-f(n))^2+2(\psi(n)-An)f(n),
\end{align*}
we now have, on summing over $n\le x$, 
\begin{align*}
\frac{Bx^3}3+o(x^3)&\le\frac{A^2x^3}3+\frac{4\epsilon^2C^2x^3}3+o(x^3). 
\end{align*}
If $\epsilon=\epsilon(A,B,C)$ is sufficiently small, and $x$ is large, then the inequality here is untenable for $A^2<B$. The theorem is proved.

\section{Segal's theorem on $\phi(n)$}
{\bf Lemma}. \emph{For Euler's function $\phi(n)$, we have, as $x\to\infty$, 
\begin{equation}\label{mertens}\sum_{n\le x}\phi(n)=\frac{Ax^2}2+O(x\log x)\end{equation}
and 
\begin{equation}\label{segal}\sum_{n\le x}\phi^2(n)=\frac{Bx^3}3+O(x^2\log^2x),\end{equation}
where, for primes $p$, 
\begin{equation*}\label{aandb}A=\prod_p\Big(1-\frac1{p^2}\Big)\qquad\text{and}\qquad
B=\prod_p\Big(1-\frac2{p^2}+\frac1{p^3}\Big). \end{equation*}}

Thus $\phi\in\mm$, and it is readily seen that $A^2<B$, so that $\phi(n)$ does not have a normal order. 
The asymptotic formula~(\ref{mertens}) is due to F.~Mertens~\cite{fm}, and (\ref{segal}) is due to~Segal~\cite{sls}, who gave a somewhat elaborate proof. 
For completeness sake, we give the proof of the lemma here. 

\medskip
\emph{Proof}. 
By M\"obius inversion, we have
$$\frac{\phi(n)}n=\sum_{d|n}\frac{\mu(d)}d,$$
where $\mu(n)$ is the M\"obius function; the formula can also be verified by taking $n$ to be a prime power, and noting that the functions involved are multiplicative. 
It follows that, as $x\to\infty$, 
\begin{align*}
\sum_{n\le x}\phi(n)&=\sum_{ab\le x}a\mu(b)=\sum_{b\le x}\mu(b)\sum_{a\le x/b}a
=\sum_{b\le x}\mu(b)\Big\{\frac12\Big(\frac{x}{b}\Big)^2+O\Big(\frac{x}{b}\Big)\Big\}\\
&=\frac{Ax^2}2+E_1(x)+E_2(x),
\end{align*}
where
\begin{align*}
A&=\sum_{b=1}^\infty\frac{\mu(b)}{b^2}=\prod_p\Big(1-\frac1{p^2}\Big),\\
E_1(x)&=O\Big(x^2\sum_{b>x}\frac1{b^2}\Big)=O(x),\qquad
E_2(x)=O\Big(x\sum_{b\le x}\frac1b\Big)=O(x\log x),
\end{align*}
so that~(\ref{mertens}) is proved.

Again, from the functions involved being multiplicative, it can be checked that
\begin{equation*}
\Big(\sum_{d|n}\frac{\mu(d)}{d}\Big)^2=\sum_{a|n}\frac{\mu^2(a)}{a^2}g(a),
\qquad\text{where}\quad g(a)=\prod_{p|a}(1-2p).
\end{equation*}
Thus, as $x\to\infty$, 
\begin{align*}
\sum_{n\le x}\phi^2(n)&=\sum_{ab\le x}a^2\mu^2(b)g(b)
=\sum_{b\le x}\mu^2(b)g(b)\Big\{\frac{x^3}{3b^3}+O\Big(\frac{x^2}{b^2}\Big)\Big\}\\
&=\frac{Bx^3}3+E_3(x)+E_4(x)
\end{align*}
where
\begin{align*}
B&=\sum_{b=1}^\infty\frac{\mu^2(b)g(b)}{b^3}=\prod_p\Big(1-\frac2{p^2}+\frac1{p^3}\Big),\\
E_3(x)&=O\Big(x^3\sum_{b>x}\frac{|g(b)|}{b^3}\Big),\qquad
E_4(x)=O\Big(x^2\sum_{b\le x}\frac{|g(b)|}{b^2}\Big).
\end{align*}
Apply the bound $|g(b)|\le\prod_{p|b}(2p)\le2^{\omega(b)}b\le d(b)b$, where $d(n)$ is the divisor function,  
and consider
\begin{align*}
\sum_{b>x}\frac{d(b)}{b^2}&=\sum_{uv>x}\frac1{u^2v^2}=\sum_{u\le x}\frac1{u^2}\sum_{v>x/u}\frac1{v^2}
+\sum_{u>x}\frac1{u^2}\sum_{v=1}^\infty\frac1{v^2}\\
&=O\Big(\frac1x\sum_{u\le x}\frac1u\Big)+O\Big(\sum_{u>x}\frac1{u^2}\Big)=O\Big(\frac{\log x}x\Big),\\
\sum_{b\le x}\frac{d(b)}{b}&=\sum_{uv\le x}\frac1{uv}=O(\log^2x).
\end{align*}
Thus $E_3(x)=O(x^2\log x)$ and $E_4(x)=O(x^2\log^2x)$, and the lemma is proved. 

\medskip
Finally, we remark that Tur\'an's method is more flexible than what is used to establish the theorem. 
Roughly speaking, the argument applies to any $\psi(n)$ for which the second moment sum $\sum_{n\le x}\psi^2(n)$
is substantially larger than what `might be expected' from the bound for the first moment sum $\sum_{n\le x}\psi(n)$. 
For example, from  
\begin{equation*}
\sum_{n\le x}d(n)\sim x\log x\qquad\text{and}\qquad\sum_{n\le x}d^2(n)\sim\frac{x\log^3x}{\pi^2},
\qquad\text{as}\quad x\to\infty,\end{equation*}
we see that the average value for $d(n)$ is~$\log n$, whereas the average value for $d^2(n)$ is~$\log^3n/\pi^2$, which is significantly larger than~$\log^2n$. 
The proof of the theorem can easily be adapted to show that $d(n)$ does not have a normal order.

\medskip


\end{document}